\documentclass{amsart}
\usepackage{amsmath,amssymb}
     \newtheorem{theorem}{Theorem}[section]
     \newtheorem{remark}[theorem]{Remark}
     \newtheorem{proposition}[theorem]{Proposition}

     \newcommand{\ct}[1]{\langle {#1}\rangle \lower.3ex\hbox{$_{t}$}}
     \newcommand{\lt}[1]{[ {#1}] \lower.3ex\hbox{$_{t}$}}

\begin{document}

\title[The $p$-Faber-Krahn Inequality Noted]{The $p$-Faber-Krahn Inequality Noted}





\author{Jie Xiao}
\address{Department of Mathematics and Statistics, Memorial University of Newfoundland, St. John's, NL A1C 5S7, Canada}
\email{jxiao@mun.ca}


\date{}



\begin{abstract}
When revisiting the Faber-Krahn inequality for the
principal $p$-Laplacian eigenvalue of a bounded open set in $\mathbb
R^n$ with smooth boundary, we simply rename it as the
$p$-Faber-Krahn inequality and interestingly find that this
inequality may be improved but also characterized through
Maz'ya's capacity method, the Euclidean volume, the Sobolev type
inequality and Moser-Trudinger's inequality.
\end{abstract}
\maketitle

\section{The $p$-Faber-Krahn Inequality Introduced}

Throughout this article, we always assume that $\Omega$ is a bounded
open set with smooth boundary $\partial\Omega$ in the $2\le
n$-dimensional Euclidean space $\mathbb R^n$ equipped with the
scalar product $\langle\cdot,\cdot\rangle$, but also $dV$ and $dA$
stand respectively for the $n$ and $n-1$ dimensional Hausdorff
measure elements on $\mathbb R^n$. For $1\le p<\infty$, the
$p$-Laplacian of a function $f$ on $\Omega$ is defined by
$$
\Delta_{p} f=-\operatorname{div}\,(|\nabla f|^{p-2}\nabla f).
$$
As usual, $\nabla$ and $
\operatorname{div}\,(|\nabla|^{p-2}\nabla)$ mean the
gradient and $p$-harmonic operators respectively (cf. \cite{DoI}).
If $W_0^{1,p}(\Omega)$ denotes the $p$-Sobolev space on $\Omega$ --
the closure of all smooth functions $f$ with compact support in
$\Omega$ (written as $f\in C^\infty_0(\Omega))$ under the norm
$$
\Big(\int_\Omega|f|^pdV\Big)^{1/p}+\Big(\int_\Omega|\nabla
f|^pdV\Big)^{1/p},
$$
then the principal $p$-Laplacian eigenvalue of $\Omega$ is defined
by
$$
\lambda_p(\Omega):=\inf\left\{\frac{\displaystyle{\int_\Omega |\nabla f|^pdV}}
{\displaystyle{\int_\Omega |f|^pdV}}:\ 0\not=f\in W_0^{1,p}(\Omega)\right\}.
$$
This definition is justified by the well-known fact that
$\lambda_2(\Omega)$ is the principal eigenvalue of the positive
Laplace operator $\Delta_2$ on $\Omega$ but also two kinds of
observation that are made below. One is the normal setting: If $p\in
(1,\infty)$, then according to \cite{Sak} there exists a nonnegative
function $u\in W^{1,p}_0(\Omega)$ such that the Euler-Lagrange
equation
\vskip5pt
\centerline{$
\Delta_p u-\lambda_p(\Omega)|u|^{p-2}u=0~~\text{in}~ \Omega$}
\vskip5pt
\noindent
holds in the weak sense of
$$
\int_\Omega \langle |\nabla u|^{p-2}\nabla u,\nabla\phi\rangle
dV=\lambda_p(\Omega)\int_{\Omega}|u|^{p-2}u\phi dV\quad\forall\
\phi\in C^\infty_0(\Omega).
$$
The other is the endpoint setting: If $p=1$, then since
$\lambda_1(\Omega)$ may be also evaluated by
$$
\inf\left\{\frac{\displaystyle{\int_\Omega |\nabla
f|dV+\int_{\partial\Omega}|f|dA}}{\displaystyle{\int_\Omega |f|dV}}:
 0\not=f\in
BV(\Omega)\right\},
$$
where $BV(\Omega)$, containing $W^{1,1}_0(\Omega)$, stands for the
space of functions with bounded variation on $\Omega$ (cf.
\cite[Chapter 5]{EvGa}), according to \cite[Theorem 4]{Dem} (cf.
\cite{HebSai}) there is a nonnegative function $u\in BV(\Omega)$
such that
\vskip5pt
\centerline{$
\Delta_1 u-\lambda_1(\Omega)|u|^{-1}u=0~~\text{in}~~ \Omega$}
\vskip5pt
\noindent
in the sense that there exists a vector-valued function
$\sigma:\Omega\mapsto \mathbb R^n$ with
\vskip5pt
\centerline{$
\|\sigma\|_{L^\infty(\Omega)}=\inf\{c:\ |\sigma|\le c~~\text{a.e.
in}~ \Omega\}<\infty$}
\vskip5pt
\noindent
and
\begin{align*}
{}& \operatorname{div}\,(\sigma)=\lambda_{1}(\Omega),
\\
&\langle\sigma,\nabla u\rangle=|\nabla u|~~\text{in}~~ \Omega,
\\
&\langle\sigma,{\mathbf{n}}\rangle u=-|u|~~\text{on}~~ \partial\Omega,
\end{align*}
where $\bf n$ represents the unit outer normal vector along
$\partial\Omega$. Moreover, it is worth pointing out that
\begin{equation}\label{eqL}
\lambda_1(\Omega)=\lim_{p\to\infty}\lambda_p(\Omega),
\end{equation}
and so that $\Delta_1 u=\lambda_1(\Omega)|u|^{-1}u$ has no classical
nonnegative solution in $\Omega$: In fact, if not, referring to
\cite[Remark 7]{KaFr} we have that for $p>1$ and $|\nabla u(x)|>0$
\begin{align}
\Delta_p u(x)
&
=(1-p)|\nabla u(x)|^{p-4}\langle D^2 u(x)\nabla
u(x),\nabla u(x)\rangle
\notag\\
&
+(n-1)H(x)|\nabla u(x)|^{p-1},
\label{eq0a}
\end{align}
where $D^2u(x)$ and $H(x)$ are the Hessian matrix of $u$ and the
mean curvature of the level surface of $u$ respectively, whence
getting by letting $p\to 1$ in (\ref{eq0a}) that
$(n-1)H(x)=\lambda_1(\Omega)$ -- namely all level surfaces of $u$
have the same mean curvature $\lambda_1(\Omega)(n-1)^{-1}$ -- but
this is impossible since the level sets $\{x\in\Omega: u(x)\ge t\}$
are strictly nested downward with respect to $t>0$.

Interestingly, Maz'ya's \cite[Theorem 8.5]{Maz0a} tells us that
$\lambda_p(\Omega)$ has an equivalent description below:
\begin{equation}\label{eq10}
\lambda_p(\Omega)\le\gamma_p(\Omega):=\inf_{\Sigma\in
AC(\Omega)}{\operatorname{cap}_p(\bar{\Sigma};\Omega)}{V(\Sigma)^{-1}}\le{p^p}{(p-1)^{1-p}}\lambda_p(\Omega).
\end{equation}
Here and henceforth, for an open set $O\subseteq\mathbb R^n$,
$AC(O)$ stands for the admissible class of all open sets $\Sigma$
with smooth boundary $\partial\Sigma$ and compact closure
$\bar{\Sigma}\subset\Omega$, and moreover
$$
\operatorname{cap}_p(K;O):=\inf\Big\{\int_O|\nabla f(x)|^pdx:\ f\in
C_0^\infty(O)\quad\&\quad f\ge 1\quad\text{in}~~ K\Big\}
$$
represents the $p$-capacity of a compact set $K\subset O$ relative
to $O$ -- this definition is extendable to any subset $E$ of $O$ via
\vskip5pt
\centerline{$
\operatorname{cap}_p(E;O):=\sup\{\operatorname{cap}_p(K;O):
\text{compact}~K\subseteq E\}$}
\vskip5pt
\noindent
-- of particular interest is that a combination of Maz'ya's \cite[p.
107, Lemma]{Maz0} and the H\"older inequality yields
\begin{equation}\label{eqlimit}
\operatorname{cap}_1(E;O)=\lim_{p\to 1}\operatorname{cap}_p(E;O).
\end{equation}
The constant $\gamma_p(\Omega)$ is called the $p$-Maz'ya constant of
$\Omega$. Of course, if $p=1$, then $(p-1)^{p-1}$ is taken as $1$ and
hence the equalities in (\ref{eq10}) are valid -- this situation
actually has another description (cf. Maz'ya \cite{Maz1}):
\begin{equation}\label{eq101}
\lambda_1(\Omega)=\gamma_1(\Omega)=h(\Omega):=\inf_{\Sigma\in
AC(\Omega)}{A(\partial\Sigma)}{V(\Sigma)^{-1}}.
\end{equation}
The right-hand-side constant in (\ref{eq101}) is regarded as the
Cheeger constant of $\Omega$ which has a root in \cite{Che}. As an
extension of Cheeger's theorem in \cite{Che}, Lefton and Wei
\cite{LeWe} (cf. \cite{KaFr} and \cite{Gri}) obtained the following
inequality:
\begin{equation}\label{eq11}
\lambda_p(\Omega)\ge p^{-p}h(\Omega)^p.
\end{equation}
Generally speaking, the reversed inequality of (\ref{eq11}) is not
true at all for $p>1$. In fact, referring to Maz'ya's first example
in \cite{Maz1}, we choose $Q$ to be the open $n$-dimensional unit
cube centered at the origin of $\mathbb R^n$. If $K$ is a compact
subset of $Q$ with $A(K)=0$ and $\operatorname{cap}_p(K;\mathbb R^n)>0$, and
if $\Omega=\mathbb R^n\setminus \cup_{z\in\mathbb Z^n}(K+z)$, i.e.,
the complement of the union of all integer shifts of $K$, then
$h(\Omega)=\gamma_1(\Omega)=0$ and $\lambda_p(\Omega)>0$ thanks to
Maz'ya's \cite[p.425, Theorem]{Maz0}, and hence there is no constant
$c_1(p,n)>0$ only depending on $1<p<n$ such that
$\lambda_p(\Omega)\le c_1(p,n) h(\Omega)^p$. Moreover, Maz'ya's
second example in \cite{Maz1} shows that if $\Omega$ is a subdomain
of the unit open ball $B_1(o)$ of $\mathbb R^n$, star-shaped with
respect to an open ball $B_\rho(o)\subset\mathbb R^n$ centered at
the origin $o$ with radius $\rho\in (0,1)$, then there is no constant
$c_2(p,n)>0$ depending only on $1<p\le n-1$ such that
$\lambda_p(\Omega)\le c_2(p,n) h(\Omega)^p$.

Determining the principal $p$-Laplacian eigenvalue of $\Omega$ is,
in general, a really hard task that relies on the value of $p$ and
the geometry of $\Omega$. However, the Faber-Krahn inequality for
this eigenvalue of $\Omega$, simply called the $p$-Faber-Krahn
inequality, provides a good way to carry out the task. To be more
precise, let us recall the content of the $p$-Faber-Krahn
inequality: If $\Omega^\ast$ is the Euclidean ball with the same
volume as $\Omega$'s, i.e., $V(\Omega^\ast)=V(\Omega)=r^n\omega_n$
(where $\omega_n$ is the volume of the unit ball in $\mathbb R^n$),
then
\begin{equation}\label{eq0b}
\lambda_p(\Omega)\ge\lambda_p(\Omega^\ast)
\end{equation}
for which equality holds if and only if $\Omega$ is a ball. A proof
of (\ref{eq0b}) can be directly obtained by Schwarz's symmetrization
-- see for example \cite[Theorem 1]{KaFr}, but the equality
treatment is not trivial -- see \cite{Bha} for an argument. Of
course, the case $p=2$ of this result goes back to the well-known
Faber-Krahn inequality (see also \cite[Theorem III.3.1]{Cha} for
an account) with $\lambda_2(\Omega^\ast)$ being $(j_{(n-2)/2}/r)^2$,
where $j_{(n-2)/2}$ is the first positive root of the Bessel
function $J_{(n-2)/2}$ and $r$ is the radius of $\Omega^\ast$. Very
recently, in \cite{Maz1} Maz'ya used his capacitary techniques to
improve the foregoing special inequality. Such a paper of Maz'ya and
his other two \cite{Maz0a}--\cite{Maz0b}, together with some
Sobolev type inequalities for
$\lambda_2(\Omega)\ge\lambda_2(\Omega^\ast)$ described in
\cite[Chapter VI]{Cha}, motivate our consideration of not only a
possible extension of Maz'ya's result -- for details see Section 2
of this article, but also some interesting geometric-analytic
properties of (\ref{eq0b}) -- for details see Section 3 of this
article.

\section{The $p$-Faber-Krahn Inequality Improved}

In order to establish a version stronger than (\ref{eq0b}), let us
recall that if from now on $B_r(x)$ represents the Euclidean ball
centered at $x\in\mathbb R^n$ of radius $r>0$, then (cf.
\cite[p.~106]{Maz0})
\begin{equation}\label{eq0}
\hskip-2pt
\operatorname{cap}_p({B_r(x)};O)\!=\!\!
\begin{cases}
n\omega_n\big(\frac{n-p}{p-1}\big)^{p-1}r^{n-p}&
\text{when}~~
O=\mathbb R^n\;\&\;p\in [1,n),
\\
0
&
\text{when}~
O=B_r(x)\;\&\;p=n,
\\
n\omega_n\big(\frac{p-n}{p-1}\big)^{p-1}r^{n-p}
&
\text{when}~
O=B_r(x)\;\&\; p\in (n,\infty).
\end{cases}
\end{equation}

\begin{proposition}\label{th1} For $t\in (0,\infty)$ and
$f\in C_0^\infty(\Omega)$, let
$\Omega_t=\{x\in\Omega: |f(x)|\ge t\}$.
\vskip5pt

\item{\rm(i)} If $p=1$, then
$$
\lambda_1(\Omega^\ast)\le\frac{\displaystyle{(n\omega_n^\frac1n)^\frac{n}{n-1}\int_\Omega|\nabla
f|dV}}{\displaystyle{\int_0^\infty\min\{\operatorname{cap}_1\big(\Omega^\ast;\mathbb
R^n\big)^\frac{n}{n-1},\operatorname{cap}_1\big(\Omega_t;\Omega\big)^\frac{n}{n-1}\}
dt}}.
$$

\item{\rm(ii)} If $p\in (1,n)$, then
$$
\lambda_p(\Omega^\ast)\le\frac{\displaystyle{(n^n\omega_n^p)^\frac{1}{n-p}
\Big(\frac{n-p}{p-1}\Big)^\frac{n(p-1)}{n-p}
\int_\Omega|\nabla
f|^pdV}}{\displaystyle{\int_0^\infty{\big(\operatorname{cap}_p\big(\Omega^\ast;\mathbb
R^n\big)^\frac{1}{1-p}+\operatorname{cap}_p\big(\Omega_t;\Omega\big)^\frac1{1-p}\big)^\frac{n(1-p)}{n-p}}\,dt^p}}.
$$

\item{\rm(iii)} If $p=n$, then
$$
\lambda_n(\Omega^\ast)\le\frac{\displaystyle{V(\Omega^\ast)^{-1}\int_\Omega|\nabla
f|^ndV}}{\displaystyle{\int_0^\infty\exp\big(-n^\frac{n}{n-1}\omega_n^\frac1{n-1}\operatorname{cap}_n\big(\Omega_t;\Omega\big)^\frac{1}{1-n}\big)\,dt^n}}.
$$

\item{\rm(iv)} If $p\in (n,\infty)$, then
$$
\lambda_p(\Omega^\ast)\le\frac{
\displaystyle{(n^n\omega_n^p)^\frac1{n-p}\Big(\frac{p-n}{p-1}\Big)^\frac{n(p-1)}{n-p}\int_\Omega|\nabla
f|^pdV}}{\displaystyle{\int_0^\infty{\big(\operatorname{cap}_p\big(\Omega^\ast;\Omega^\ast\big)^\frac{1}{1-p}-\operatorname{cap}_p\big(\Omega_t;\Omega\big)^\frac{1}{1-p}\big)^\frac{n(1-p)}{p-n}}\,dt^p}}.
$$

\item{\rm(v)} The inequalities in {\rm(i)}--{\rm(ii)}--{\rm(iii)}--{\rm(iv)} imply the inequality {\rm(\ref{eq0b})}.
\end{proposition}

\begin{proof} For simplicity, suppose that $r=(V(\Omega)\omega_n^{-1})^\frac1n$ is the radius of the Euclidean ball $\Omega^\ast$,
$\Omega_t^\ast$ is the Euclidean ball with
$V(\Omega_t^\ast)=V(\Omega_t)$, and $f^\ast$ equals
$\displaystyle{\int_0^\infty
1_{\Omega_t^\ast}dt}$, where $1_{E}$ stands for the characteristic
function of a set $E\subseteq\mathbb R^n$. Then
$$
\int_\Omega|\nabla f^\ast|^pdV\le\int_\Omega|\nabla
f|^pdV\quad\hbox{\&}\quad\int_\Omega |f^\ast|^pdV=\int_\Omega
|f|^pdV.
$$
Consequently, from the definitions of $\lambda_p(\Omega^\ast)$ and
$f^\ast$ as well as \cite[p.38, Exercise 1.4.1]{Dac} it follows that
\begin{equation}\label{eq1}
\lambda_p(\Omega^\ast)\int_0^r|a(t)|^p t^{n-1}dt\le\int_0^r|a'(t)|^p
t^{n-1}dt
\end{equation}
holds for any absolutely continuous function $a$ on $(0,r]$ with
$a(r)=0$.

{\it Case 1.} Under $p\in (1,n)$, set
$$
s=\frac{t^\frac{p-n}{p-1}-r^\frac{p-n}{p-1}}{\alpha},\quad\text{where}~~
\alpha=(n\omega_n)^\frac{1}{p-1}\Big(\frac{n-p}{p-1}\Big).
$$
This yields
$$
t=(r^\frac{p-n}{p-1}+\alpha s)^\frac{p-1}{p-n}\quad\text{and}\quad
\frac{dt}{ds}=\frac{\alpha(p-1)}{p-n}\Big(\alpha
s+r^\frac{p-n}{p-1}\Big)^\frac{n-1}{p-n}.
$$
If $b(s)=a(t)$, then
$$
\int_0^r|a(t)|^p
t^{n-1}dt=\Big(\frac{\alpha(p-1)}{n-p}\Big)\int_0^\infty|b(s)|^p\big(r^\frac{p-n}{p-1}+\alpha
s\big)^\frac{p(n-1)}{p-n}ds
$$
and
$$
\int_0^r|a'(t)|^pt^{n-1}dt=\Big(\frac{\alpha(p-1)}{n-p}\Big)^{1-p}\int_0^\infty|b'(s)|^pds.
$$
Consequently, (\ref{eq1}) amounts to
\begin{equation}\label{eq2}
\lambda_p(\Omega^\ast)\Big(\frac{\alpha(p-1)}{n-p}\Big)^{p}\int_0^\infty|b(s)|^p\big(r^\frac{p-n}{p-1}+\alpha
s\big)^\frac{p(n-1)}{p-n}ds\le\int_0^\infty|b'(s)|^pds.
\end{equation}

{\it Case 2.} Under $p=n$, set
$$
s=\frac{\ln\frac{r}{t}}{\beta},\quad\text{where}~~\beta=(n\omega_n)^\frac1{n-1}.
$$
This gives
$$
t=r\exp(-\beta s)\quad\text{and}\quad \frac{dt}{ds}=-\beta
r\exp(-\beta s).
$$
If $b(s)=a(t)$, then
$$
\int_0^r|a(t)|^n t^{n-1}dt=\beta
r^n\int_0^\infty|b(s)|^n\exp(-n\beta s)ds
$$
and
$$
\int_0^r|a'(t)|^nt^{n-1}dt=\beta^{1-n}\int_0^\infty|b'(s)|^nds.
$$
As a result, (\ref{eq1}) is equivalent to
\begin{equation}\label{eq3}
\lambda_n(\Omega^\ast)\beta^{n}\int_0^\infty|b(s)|^n\exp(-n\beta
s)ds\le\int_0^\infty|b'(s)|^nds.
\end{equation}

{\it Case 3.} Under $p\in (n,\infty)$, set
$$
s=\frac{r^\frac{p-n}{p-1}-t^\frac{p-n}{p-1}}{\gamma},\quad\text{where}~~
\gamma=(n\omega_n)^\frac{1}{p-1}\Big(\frac{p-n}{p-1}\Big).
$$
This produces
$$
t=(r^\frac{p-n}{p-1}-\gamma
s)^\frac{p-1}{p-n}\quad\hbox{and}\quad\frac{dt}{ds}=\Big(\frac{\gamma(p-1)}{n-p}\Big)(r^\frac{p-n}{p-1}-\gamma
s)^\frac{n-1}{p-n}.
$$
If $b(s)=a(t)$, then
$$
\int_0^r|a(t)|^p
t^{n-1}dt=\Big(\frac{\gamma(p-1)}{p-n}\Big)\int_0^{\frac{r^\frac{p-n}{p-1}}{\gamma}}|b(s)|^p
\big(r^\frac{p-n}{p-1}-\gamma s\big)^\frac{p(n-1)}{p-n}ds
$$
and
$$
\int_0^r|a'(t)|^pt^{n-1}dt=\Big(\frac{\gamma(p-1)}{p-n}\Big)^{1-p}\int_0^{\frac{r^\frac{p-n}{p-1}}{\gamma}}
|b'(s)|^pds.
$$
Thus, (\ref{eq1}) can be reformulated as

\begin{equation}\label{eq4}
\lambda_p(\Omega^\ast)\Big(\frac{\gamma(p-1)}{p-n}\Big)^{p}\int_0^{\frac{r^\frac{p-n}{p-1}}{\gamma}}
|b(s)|^p\big(r^\frac{p-n}{p-1}-\gamma
s\big)^\frac{p(n-1)}{p-n}ds\le\int_0^{\frac{r^\frac{p-n}{p-1}}{\gamma}}|b'(s)|^pds.
\end{equation}

In the three inequalities (\ref{eq2})-(\ref{eq3})-(\ref{eq4}),
choosing
$$
s=\int_0^\tau \Big(\int_{\{x\in\Omega: f(x)=t\}}|\nabla
f|^{p-1}dA\Big)^\frac{1}{1-p}dt
$$
and letting $\tau(s)$ be the inverse of the last function, we have
two equalities:
\begin{equation}\label{eq5}
\frac{ds}{d\tau}=\frac1{\tau'(s)}\quad\&\quad \int_0^\infty
|s'(\tau)|^{-p}d\tau=\int_\Omega|\nabla f|^pdV
\end{equation}
and Maz'ya's inequality for the $p$-capacity (cf.
\cite[p.102]{Maz0}):
\begin{equation}\label{eq6}
s\le\operatorname{cap}_p\big(\Omega_{\tau(s)};\Omega\big)^\frac1{1-p}.
\end{equation}
The above estimates (\ref{eq0}) and
(\ref{eq2})-(\ref{eq3})-(\ref{eq4})-(\ref{eq5})-(\ref{eq6}) give the
inequalities in (ii)-(iii)-(iv).

Next, we verify (i). In fact, this assertion follows from
formulas (\ref{eqL}) and (\ref{eqlimit}), taking the limit $p\to 1$
in the inequality established in (ii), and using the elementary
limit evaluation
$$
\lim_{p\to
1}(c_1^\frac1{p-1}+c_2^\frac1{p-1})^{p-1}=\max\{c_1,c_2\}\quad\hbox{for}\quad
c_1,c_2\ge 0.
$$

Finally, we show (v). To do so, recall Maz'ya's lower bound
inequality for $\operatorname{cap}_p(\cdot,\cdot)$ (cf. \cite[p.105]{Maz0}):
\begin{equation}\label{eq7}
\operatorname{cap}_p(\Omega_t;\Omega)\ge\Big(\int_{V(\Omega_t)}^{V(\Omega)}\mu(v)^\frac{p}{1-p}dv\Big)^{1-p}\quad\hbox{for}\quad
0<t, p-1<\infty,
\end{equation}
where $\mu(v)$ is defined as the infimum of $A(\partial\Sigma)$ over
all open subsets $\Sigma\in AC(\Omega)$ with $V(\Sigma)\ge v$.

From the classical isoperimetric inequality with sharp constant
\begin{equation}\label{eq7f}
V(\Sigma)^\frac{n-1}{n}\le(n\omega_n^\frac{1}{n})^{-1}A(\partial\Sigma)\quad\forall\
\Sigma\in AC(\mathbb R^n)
\end{equation}
it follows that $\mu(v)\ge n\omega_n^\frac1n v^\frac{n-1}{n}$ and
consequently
\begin{equation}\label{eq8}
\int_{V(\Omega_t)}^{V(\Omega)}\mu(v)^\frac{p}{1-p}dv
\le
\begin{cases}
\displaystyle{\frac{V(\Omega_t)^\frac{p-n}{n(p-1)}-V(\Omega)^\frac{p-n}{n(p-1)}}
{\Big(\frac{n(p-1)}{(n-p)(n\omega^{1/n})^{{p}/({p-1})}}\Big)^{-1}}}
&\text{for}~~
1<p\not=n,
\\[23pt]
\displaystyle{(n\omega_n^{1/n})^{n/(1-n)}\ln
\Big(\frac{V(\Omega)}{V(\Omega_t)}\Big)}
&\text{for}~~
p=n.
\end{cases}
\end{equation}
Using (\ref{eq8}) and (ii)-(iii) we derive the following estimates.

{\it Case 1.} If $1<p<n$, then
\begin{align*}
I_{1<p<n}&:=\int_0^\infty\Big(\operatorname{cap}_p\big(\Omega^\ast;\mathbb
R^n\big)^\frac{1}{1-p}+\operatorname{cap}_p\big(\Omega_t;\Omega\big)^\frac1{1-p}\Big)^\frac{n(p-1)}{p-n}\,dt^p
\\
&\ge\int_0^\infty\Biggl(\operatorname{cap}_p\big(\Omega^\ast;\mathbb
R^n\big)^\frac{1}{1-p}+\frac{V(\Omega_t)
^\frac{p-n}{n(p-1)}-V(\Omega)^\frac{p-n}{n(p-1)}}
{\Big(\frac{n(p-1)}{(n-p)(n\omega^{1/n})^{p/(p-1)}}\Big)^{-1}}\Biggr)^
\frac{n(p-1)}{p-n}\,dt^p
\\
&=\left(\frac{n(p-1)}{(n-p)(n\omega^\frac1n)^\frac{p}{p-1}}\right)^\frac{p-n}{n(p-1)}\int_0^\infty V(\Omega_t)\,dt^p\\
&=\left(\frac{n(p-1)}{(n-p)(n\omega^\frac1n)^\frac{p}{p-1}}\right)^\frac{p-n}{n(p-1)}\int_\Omega|f|^pdV.
\end{align*}

{\it Case 2.} If $p=n$, then
\begin{align*}
I_{p=n}&:=\int_0^\infty\exp\Big(-n^\frac{n}{n-1}\omega_n^\frac1{n-1}\operatorname{cap}_n\big(\Omega_t;\Omega\big)^\frac{1}{1-n}\Big)\,dt^n\\
&\ge
V(\Omega)^{-1}\int_0^\infty{V(\Omega_t)}\,dt^n\\
&=V(\Omega)^{-1}\int_\Omega|f|^ndV.
\end{align*}

{\it Case 3.} If $n<p<\infty$, then
\begin{align*}
I_{n<p<\infty}&:=\int_0^\infty\Big(\operatorname{cap}_p\big(\Omega^\ast;\Omega^\ast\big)^\frac{1}{1-p}-\operatorname{cap}_p\big(\Omega_t;\Omega\big)^\frac{1}{1-p}\Big)^\frac{n(p-1)}{p-n}\,dt^p\\
&\ge\int_0^{\infty}\left(\operatorname{cap}_p\big(\Omega^\ast;\Omega^\ast
\big)^\frac{1}{1-p}-\frac{V(\Omega)^\frac{p-n}{n(p-1)}
-V(\Omega_t)^\frac{p-n}{n(p-1)}}{\Big(\frac{n(p-1)}{(p-n)(n\omega_
n^{1/n})^{{p}/({p-1})}}\Big)^{-1}}\right)^\frac{(p-1)n}{p-n}\,dt^p
\\
&= \left(\frac{n(p-1)}{(p-n)(n\omega_n^\frac1n)^\frac{p}{p-1}}\right)^\frac{(p-1)n}{p-n} \int_0^{\infty}V(\Omega_t)dt^p\\
&=\left(\frac{n(p-1)}{(p-n)(n\omega_n^\frac1n)^\frac{p}{p-1}}\right)^\frac{(p-1)n}{p-n}\int_\Omega
|f|^pdV.
\end{align*}

Now the last three cases, along with (ii)-(iii)-(iv), yield (v) for
$1<p<\infty$. In order to handle the setting $p=1$, letting $p\to 1$
in (\ref{eq7}) we employ (\ref{eqlimit}) and
$$
\lim_{p\to 1}(1-c^\frac1{1-p})^{1-p}=1\quad\hbox{for}\quad c\ge 1
$$
to achieve the following relative isocapacitary inequality with
sharp constant:
\begin{equation}\label{eq9}
\operatorname{cap}_1\big(\Omega_t;\Omega\big)\ge n\omega_n^\frac1n
V(\Omega_t)^\frac{n-1}{n}.
\end{equation}
As a consequence of (\ref{eq9}), we find
\begin{eqnarray*}
I_{p=1}&:=&\int_0^\infty\min\{\operatorname{cap}_1\big(\Omega^\ast;\mathbb
R^n\big)^\frac{n}{n-1},\operatorname{cap}_1\big(\Omega_t;\Omega\big)^\frac{n}{n-1}\}\,
dt\\
&\ge&(n\omega_n^\frac1n)^\frac{n}{n-1}\int_0^\infty\min\{V(\Omega),V(\Omega_t)\}\,dt\\
&=&(n\omega_n^\frac1n)^\frac{n}{n-1}\int_0^\infty V(\Omega_t)\,dt\\
&=&(n\omega_n^\frac1n)^\frac{n}{n-1}\int_\Omega|f|dV,
\end{eqnarray*}
thereby getting the validity of (v) for $p=1$ thanks to (i).
\end{proof}

\begin{remark}\label{rm1}
Perhaps it is appropriate to mention that (ii)-(iii)-(iv) in
Proposition \ref{th1} can be also obtained through choosing $q=p\in
(1,\infty)$ and letting $\mathbf M(\theta)$-function in Maz'ya's
\cite[Theorem 2]{Maz1} be respectively
\vskip5pt


\begin{equation*}
\begin{cases}

&\displaystyle{\frac{\lambda_p(\Omega^\ast)(n^n\omega_n^p)^\frac{1}{p-n}\Big(\frac{n-p}{p-1}\Big)^\frac{n(p-1)}{p-n}}
{\Big(\operatorname{cap}_p\big(\Omega^\ast;\mathbb
R^n\big)^\frac{1}{1-p}+\theta\Big)^\frac{n(1-p)}{p-n}}\quad\hbox{for}\ p\in (1,n),}\\
&
\displaystyle{\lambda_n(\Omega^\ast)V(\Omega^\ast)
\exp\big(-(n^n\omega_n)^\frac1{n-1}\theta\big)
~\hbox{for}~p=n},
\\
&
\displaystyle{\frac{\lambda_p(\Omega^\ast)(n^n\omega_n^p)^\frac{1}{p-n}\Big(\frac{p-n}{p-1}\Big)^\frac{n(p-1)}{p-n}}
{\big(\operatorname{cap}_p\big(\Omega^\ast;\mathbb
R^n\big)^\frac{1}{1-p}-\theta\big)^\frac{n(1-p)}{p-n}}
\hskip7mm
\hbox{for}~
\theta\le\operatorname{cap}_p\big(\Omega^\ast;\mathbb R^n\big)^\frac{1}{1-p}
~\&~p\in (n,\infty)},
\\
&
\displaystyle{0\hskip46mm \hbox{for}~\theta>\operatorname{cap}_p\big(\Omega^\ast;\mathbb
R^n\big)^\frac{1}{1-p}~ \&~ p\in (n,\infty).}
\end{cases}
\end{equation*}
\end{remark}

\section{The $p$-Faber-Krahn Inequality Characterized}

When looking over the $p$-Faber-Krahn inequality (\ref{eq0b}), we
get immediately its alternative (cf. \cite{FuMaPr1,FuMaPr2})
as follows:
\begin{equation}\label{eq12}
\lambda_p(\Omega){V(\Omega)}^{\frac{p}{n}}\ge{\lambda_p\big(B_1(o)\big)}{\omega_n^\frac{p}{n}}.
\end{equation}
It is well known that (\ref{eq12}) is sharp in the sense that if
$\Omega$ is a Euclidean ball in $\mathbb R^n$, then equality of
(\ref{eq12}) is valid. Although the explicit value of
$\lambda_p\big(B_1(o)\big)$ is so far unknown except
\begin{equation}\label{eq12a}
\lambda_1\big(B_1(o)\big)=n\quad\&\quad\lambda_2
\big(B_1(o)\big)=j_{(n-2)/2}^2,
\end{equation}
Bhattacharya's \cite[Lemma 3.4]{Bha}
yields
\begin{equation}\label{eq13}
\lambda_p\big(B_1(o)\big)\ge n^{2-p}{p}^{p-1}(p-1)^{1-p},
\end{equation}
whence giving $\lambda_1\big(B_1(o)\big)\ge n$. Meanwhile, from
Proposition \ref{th1} we can get an explicit upper bound of
$\lambda_p\big(B_1(o)\big)$ via selecting a typical test function in
$W^{1,p}_0\big(B_1(o)\big)$, particularly finding
$\lambda_1\big(B_1(o)\big)\le n$ and hence the first formula in
(\ref{eq12a}).

Although it is not clear whether Colesanti--Cuoghi--Salani's geometric
Brunn--Minkowski type inequality of $\lambda_p(\Omega)$ for convex
bodies $\Omega$ in \cite{CoCuSa} can produce (\ref{eq12}), a
geometrical-analytic look at (\ref{eq12}) leads to the forthcoming
investigation in accordance with four situations: $p=1$; $1<p<n$;
$p=n$; $n<p<\infty$.

The case $p=1$ is so special that it produces sharp geometric and
analytic isoperimetric inequalities indicated below.

\begin{proposition}\label{th2} The following statements are
equivalent{\rm:}

\item{\rm(i)} The sharp $1$-Faber-Krahn inequality
$$
\lambda_1(\Omega)V(\Omega)^\frac1n\ge n\omega_n^\frac1n\quad\forall\
\Omega\in AC(\mathbb R^n)
$$
holds.

\item{\rm(ii)} The sharp $(1,\frac{1-n}{n})$-Maz'ya isocapacitary
inequality
$$
\operatorname{cap}_1(\bar{\Omega};\mathbb R^n)V(\Omega)^\frac{1-n}{n}\ge
n\omega_n^\frac1n\quad\forall\ \Omega\in AC(\mathbb R^n)
$$
holds.

\item{\rm(iii)} The sharp $(1,\frac{n}{n-1})$-Sobolev inequality
$$
\Big({\int_{\mathbb R^n}|\nabla f|dV}\Big){\Big(\int_{\mathbb
R^n}|f|^\frac{n}{n-1}dV\Big)^\frac{1-n}{n}}\ge
n\omega_n^\frac1n\quad\forall\ f\in C_0^\infty(\mathbb R^n)
$$
holds.
\end{proposition}

\begin{proof} (i)$\Rightarrow$(ii) Noticing
$$
V(\Omega)^{-1}{A(\partial\Omega)}\ge\lambda_1(\Omega)\quad\forall\
\Omega\in AC(\mathbb R^n),
$$
we get (i)$\Rightarrow$(\ref{eq7f}). By Maz'ya's formula in \cite[p.
107, Lemma]{Maz0} saying
$$
\operatorname{cap}_1(\bar{\Omega};\mathbb
R^n)=\inf_{\bar{\Omega}\subset\Sigma\in AC(\mathbb
R^n)}A(\partial\Sigma)\quad\forall\ \Omega\in AC(\mathbb R^n),
$$
we further find (\ref{eq7f})$\Rightarrow$(ii).

(ii)$\Rightarrow$(iii) Under (ii), we use the end-point case of Maz'ya's inequality in \cite[Proposition 1]{Maz0b} (cf. \cite[Theorems 1.1-1.2]{Xia}) to obtain
\begin{eqnarray*}
\int_{\mathbb R^n}|f|^\frac{n}{n-1}\,dV&=&\int_0^\infty V\big(\{x\in\mathbb R^n:\ |f(x)|\ge t\}\big)\,dt^\frac{n}{n-1}\\
&\le&\int_0^\infty\Big((n\omega_n^\frac1n)^{-1}\hbox{cap}_1\big(\{x\in\mathbb R^n:\ |f(x)|\ge t\}\big)\Big)^\frac{n}{n-1}\,dt^\frac{n}{n-1}\\
&\le&(n\omega_n^\frac1n)^\frac{n}{1-n}\Big(\int_{\mathbb R^n}|\nabla f|\,dV\Big)^\frac{n}{n-1},
\end{eqnarray*}
whence getting (iii).

(iii)$\Rightarrow$(i) For $\Omega\in AC(\mathbb R^n)$ and $f\in W^{1,1}_0(\Omega)$, define a Sobolev function $g$ on $\mathbb R^n$ via putting $g=f$ in $\Omega$ and $g=0$ in $\mathbb R^n\setminus\Omega$. If (iii) holds, then the inequality in (iii) is valid for $g$. Using H\"older's inequality we have
$$
\int_\Omega |f|\,dV\le \Big(\int_\Omega|f|^\frac{n}{n-1}\,dV\Big)^\frac{n-1}{n}V(\Omega)^\frac1n
$$
and consequently,
$$
\frac{\int_\Omega|\nabla f|\,dV}{\int_\Omega |f|\,dV}\ge\frac{\int_{\mathbb R^n}|\nabla g|\,dV}{\big(\int_{\mathbb R^n} |g|^\frac{n}{n-1}\,dV\big)^\frac{n-1}{n}V(\Omega)^\frac1n}\ge\frac{n\omega_n^\frac1n}{V(\Omega)^\frac1n}.
$$
This, along with the definition of $\lambda_1(\Omega)$, yields the inequality in (i).
\end{proof}

\begin{remark}\label{rm2} $n\omega_n^\frac1n$ is the best constant
for (i)-(ii)-(iii) whose equalities occur when $\Omega=B_1(o)$ and
$f\to 1_{B_1(o)}$. Moreover, the equivalence between the classical isoperimetric inequality (\ref{eq7f}) and the Sobolev inequality (iii) above is well known and due to Federer--Fleming \cite{FF} and
Maz'ya \cite{Maz00}.
\end{remark}

Nevertheless, the setting $1<p<n$ below does not yield optimal constants.

\begin{proposition}\label{th3} For $p\in (1,n)$, the
statement (i) follows from the mutually equivalent ones (ii) and (iii) below:

\item{\rm(i)} There is a constant $\kappa_1(p,n)>0$ depending only on $p$ and $n$ such that the $p$-Faber-Krahn inequality
$$
\lambda_p(\Omega)V(\Omega)^\frac{p}{n}\ge\kappa_1(p,n)\quad\forall\
\Omega\in AC(\mathbb R^n)
$$
holds.

\item{\rm(ii)} There is a constant $\kappa_2(p,n)>0$ depending only on $p$ and $n$ such that
the $(p,\frac{p-n}{n})$-Maz'ya isocapacitary inequality
$$
\operatorname{cap}_p(\bar{\Omega};\mathbb R^n)V(\Omega)^\frac{p-n}{n}\ge\kappa_2(p,n)\quad\forall\
\Omega\in AC(\mathbb R^n)
$$
holds.

\item{\rm(iii)} There is a constant $\kappa_2(p,n)>0$ depending only on $p$ and $n$ such that the $(p,\frac{pn}{n-p})$-Sobolev inequality
$$
{\Big(\int_{\mathbb R^n}|\nabla f|^pdV\Big)}{\Big(\int_{\mathbb
R^n}|f|^\frac{pn}{n-p}dV\Big)^\frac{p-n}{n}}\ge\kappa_3(p,n)\quad\forall\
f\in C_0^\infty(\mathbb R^n)
$$
holds.
\end{proposition}

\begin{proof} Note that (ii)$\Leftrightarrow$(iii) is a special case of Maz'ya's
\cite[Theorem 8.5]{Maz0a}. So it suffices to consider the following implications.

(ii)$\Rightarrow$(i) This can be seen from
\cite{Gri}. In fact, for $\Sigma\in AC(\Omega)$ and $\Omega\in AC(\mathbb R^n)$ one has
$$
\frac{\operatorname{cap}_p(\bar{\Sigma};\Omega)}{V(\Sigma)}\ge\frac{\operatorname{cap}_p(\bar{\Sigma};\mathbb R^n)}{V(\Sigma)}\ge\kappa_2(p,n)V(\Sigma)^{-\frac{p}{n}}\ge\kappa_2(p,n)V(\Omega)^{-\frac{p}{n}}
$$
and thus by (\ref{eq10}),
$$
\lambda_p(\Omega)V(\Omega)^{\frac{p}{n}}\ge(p-1)^{p-1}p^{-p}\kappa_2(p,n).
$$

(iii)$\Rightarrow$(i) Suppose now that (iii) is true. Since there exists a nonzero
minimizer $u\in W^{1,p}_0(\Omega)$ such that
$$
\int_\Omega|\nabla u|^{p-2}\langle \nabla
u,\nabla\phi\rangle\,dV=\lambda_p(\Omega)\int_\Omega|u|^{p-2}u\phi\,dV
$$
holds for any $\phi\in C_0^\infty(\Omega)$. Letting $\phi$ approach
$u$ in the above equation, extending $u$ from $\Omega$ to $\mathbb
R^n$ via defining $u=0$ on $\mathbb R^n\setminus\Omega$, and writing
this extension as $f$, we employ (iii) and the H\"older inequality to get
\begin{align*}
\lambda_p(\Omega)
&=\frac{\displaystyle{\int_\Omega|\nabla
u|^pdV}}{\displaystyle{\int_\Omega|u|^pdV}}
=\frac{\displaystyle{\int_{\mathbb R^n} |\nabla
f|^pdV}}{\displaystyle{\int_{\mathbb R^n}|f|^pdV}}
\\[7pt]
{}&\ge \kappa_3(p,n){\Big(\int_{\mathbb R^n} |f|^\frac{pn}{n-p}dV\Big)^\frac{n-p}{n}}{\Big(\int_{\mathbb R^n}|f|^pdV\Big)^{-1}}\\
&=\kappa_3(p,n){\Big(\int_\Omega |u|^\frac{pn}{n-p}dV\Big)^\frac{n-p}{n}}{\Big(\int_\Omega|u|^pdV\Big)^{-1}}\\
&\ge\kappa_3(p,n)V(\Omega)^{-\frac{p}{n}},
\end{align*}
whence reaching (i).
\end{proof}

\begin{remark}\label{rm3} It is worth remarking that the best values of $\kappa_1(p,n)$,
$\kappa_2(p,n)$, and $\kappa_3(p,n)$ are
$$
\lambda_p\big(B_1(o)\big)\omega_n^\frac{p}{n},\
n\omega_n^\frac{p}{n}\Big(\frac{n-p}{p-1}\Big)^{p-1},\
$$
and
$$
n\omega_n^\frac{p}{n}\Big(\frac{n-p}{p-1}\Big)^{p-1}\Big(\frac{\Gamma(\frac{n}{p})\Gamma(n+1-\frac{n}{p})}{\Gamma(n)}\Big)^\frac{p}{n}
$$
respectively. These constants tend to $n\omega_n^\frac1n$ as $p\to
1$.  In addition, from Carron's paper
\cite{Car} we see that (i) implies (ii) and (iii) under $p=2$, and consequently conjecture that this implication is also valid for $p\in (1,n)\setminus\{2\}$.
\end{remark}

Clearly, (ii) and (iii) in Proposition \ref{th3} cannot be naturally
extended to $p=n$. However, they have the forthcoming replacements.

\begin{proposition}\label{th4} For $\Omega\in AC(\mathbb R^n)$, the statement (i) follows from the mutually equivalent ones (ii) and (iii) below:

\item{\rm(i)} The $n$-Faber-Krahn type inequality
$$
\lambda_n(\Omega)V(\Omega)\ge \frac{n^n\omega_n}{(n-1)!}E_n(\Omega)^{-1}
$$
holds where
\begin{equation*}
\label{eqn1}
E_n(\Omega):=\sup_{f\in C^\infty_0(\Omega),\ \int_\Omega|\nabla
f|^ndV\le
1}V(\Omega)^{-1}\int_\Omega\exp\Big(\frac{|f|^\frac{n}{n-1}}{(n^n\omega_n)^\frac1{1-n}}\Big)dV.
\end{equation*}

\item{\rm(ii)} The $(n,0)$-capacity-volume inequality
\begin{equation*}\label{eqn}
{V(\Sigma)}{V(\Omega)^{-1}}\le\exp\Big(-\Big(\frac{n^n\omega_n}{\operatorname{cap}_n(\bar{\Sigma};\Omega)}\Big)^\frac1{n-1}\Big)\quad\forall\
\Sigma\in AC(\Omega)
\end{equation*}
holds.
\item{\rm(iii)} The Moser-Trudinger inequality
$E_n(\Omega)<\infty$ holds.
\end{proposition}
\begin{proof} (ii)$\Rightarrow$(iii) Suppose (ii) holds. For $f\in C^\infty_0(\Omega)$ and $t\ge 0$ with $\int_\Omega|\nabla f|^n\,dV\le 1$ let $\Omega_t=\{x\in\Omega: |f(x)|\ge t\}$. Then the layer-cake formula gives
\begin{eqnarray*}
\int_\Omega\exp\Big(\frac{|f|^\frac{n}{n-1}}{(n^n\omega_n)^\frac1{1-n}}\Big)dV&=&\int_0^\infty V(\Omega_t)\,d\exp\big((n^n\omega_n)^\frac1{n-1}t^\frac{n}{n-1}\big)\\
&\le& V(\Omega)\int_0^\infty\frac{d\exp\big((n^n\omega_n)^\frac1{n-1}t^\frac{n}{n-1}\big)}{
\exp\big((n^n\omega_n)^\frac1{n-1}\operatorname{cap}_n(\Omega_t;\Omega)^\frac1{1-n}\big)}
\end{eqnarray*}
where the last integral is finite by Maz'ya's \cite[Proposition 2]{Maz0b}. As a result, we find
$$
E_n(\Omega)\le\int_0^\infty\frac{d\exp\big((n^n\omega_n)^\frac1{n-1}t^\frac{n}{n-1}\big)}{
\exp\big((n^n\omega_n)^\frac1{n-1}\operatorname{cap}_n(\Omega_t;\Omega)^\frac1{1-n}\big)}<\infty,
$$
thereby reaching (iii).

(iii)$\Rightarrow$(ii) If (iii) holds, then $f\in C_0^\infty(\Omega)$,
$f\ge 1$ on $\bar{\Sigma}$ and $\Sigma\in AC(\Omega)$ imply
$$
\int_\Omega\Big|\nabla\Big(\frac{f}{\big(\int_\Omega|\nabla f|^n\,dV\big)^\frac1n}\Big)\Big|^n\,dV=1,
$$
and hence
\begin{eqnarray*}
V(\Omega)E_n(\Omega)&\ge&\int_\Omega
\exp\Big((n^n\omega_n)^\frac1{n-1}|f|^\frac{n}{n-1}\Big(\int_\Omega|\nabla f|^ndV\Big)^{\frac1{1-n}}\Big)dV\\
&\ge&
V(\Sigma)\exp\Big((n^n\omega_n)^\frac1{n-1}\Big(\int_\Omega|\nabla
f|^ndV\Big)^{\frac1{1-n}}\Big),
\end{eqnarray*}
whence giving (ii) through the definition of
$\operatorname{cap}_n(\bar{\Sigma};\Omega)$.

Of course, if either (ii) or (iii) is valid, then the elementary inequality
$$
\exp t\ge \frac{t^{n-1}}{(n-1)!}\quad\forall t\ge 0
$$
yields
$$
V(\Omega)E_n(\Omega)\ge \frac{n^n\omega_n}{(n-1)!}\left(\frac{\int_\Omega|f|^n\,dV}{\int_\Omega|\nabla f|^n\,dV}\right)\quad\forall f\in C^\infty_0(\Omega),
$$
whence giving (i) by the characterization of $\lambda_n(\Omega)$ in terms of $C_0^\infty(\Omega)$ -- see also \cite{KaLi}.
\end{proof}

\begin{remark}\label{rm4} The equality of (ii) happens when $\Omega$ and $\Sigma $ are concentric Euclidean balls
-- see also \cite[p.15]{Flu}. Moreover, the supremum defining $E_n(\Omega)$ becomes infinity when $(n^n\omega_n)^\frac1{n-1}$
is replaced by any larger constant -- see also \cite[p.97-98]{Flu}.
\end{remark}

Next, let us handle the remaining case $p\in (n,\infty)$.

\begin{proposition}\label{th5} For $p\in (n,\infty)$ and $\Omega\in AC(\mathbb R^n)$, the statement (i) follows from the mutually equivalent ones (ii) and (iii) below:

\item{\rm(i)} The $p$-Faber-Krahn inequality
$$
\lambda_p(\Omega)V(\Omega)^\frac{p}{n}\ge E_{n,p}(\Omega)
$$
holds where
$$
E_{n,p}(\Omega):=\inf_{f\in C^\infty_0(\Omega),\|f\|_{L^\infty(\Omega)}\le 1}V(\Omega)^{\frac{p-n}{n}}\int_\Omega|\nabla f|^p\,dV.
$$

\item{\rm(ii)} The $(p,\frac{p-n}{n})$-capacity-volume inequality
$$
\operatorname{cap}_p(\bar{\Sigma};\Omega)V(\Omega)^\frac{p-n}{n}\ge E_{n,p}(\Omega)\quad\forall\
\Sigma\in AC(\Omega)
$$
holds.

\item{\rm(iii)} The $(p,\infty)$-Sobolev inequality
$$
\Big(\int_\Omega|\nabla
f|^p\,dV\Big)\|f\|_{L^\infty(\Omega)}^{-p}V(\Omega)^{\frac{p-n}{n}}\ge
E_{n,p}(\Omega)\quad\forall\
f\in C_0^\infty(\Omega)
$$
holds.
\end{proposition}

\begin{proof} (iii)$\Rightarrow$(ii) Suppose (iii) is valid. If $f\in C^\infty_0(\Omega)$, $\Sigma\in AC(\Omega)$ and $f\ge 1$ on $\bar{\Sigma}$, then $\|f\|_{L^\infty(\Omega)}\ge 1$ and hence
$$
V(\Omega)^{\frac{p-n}{n}}\int_\Omega |\nabla f|^p\,dV\ge E_{n,p}(\Omega)\|f\|^p_{L^\infty(\Omega)}\ge E_{n,p}(\Omega).
$$
This, plus the definition of $\hbox{cap}_p(\bar{\Sigma};\Omega)$, yields
$$
V(\Omega)^{\frac{p-n}{n}}\hbox{cap}_p(\bar{\Sigma};\Omega)\ge E_{n,p}(\Omega).
$$
Namely, (ii) holds.

(ii)$\Rightarrow$(iii) Suppose (ii) is valid. For $q>p$ and $\Sigma\in AC(\Omega)$ we have
$$
\frac{\operatorname{cap}_p(\bar{\Sigma};\Omega)}{V(\Sigma)^\frac{q-p}{q}}\ge E_{n,p}(\Omega)
V(\Omega)^{\frac{p}{q}-\frac{p}{n}}.
$$
For $f\in C^\infty_0(\Omega)$ and $t\ge 0$ let $\Omega_t=\{x\in\Omega: |f(x)|\ge t\}$. According to the layer-cake formula and \cite[Proposition 1]{Maz0b}, we have
\begin{eqnarray*}
&&\int_{\Omega}|f|^\frac{pq}{q-p}\,dV\\
&=&\int_0^\infty V(\Omega_t)\,dt^\frac{pq}{q-p}\\
&\le&\Big(E_{n,p}(\Omega)V(\Omega)^{\frac{p}{q}-\frac{p}{n}}\Big)^{\frac{q}{p-q}}
\int_0^\infty\hbox{cap}_p(\Omega_t;\Omega)^\frac{q}{q-p}\,dt^\frac{pq}{q-p}\\
&\le&\Big(E_{n,p}(\Omega)V(\Omega)^{\frac{p}{q}-\frac{p}{n}}\Big)^{\frac{q}{p-q}}\left(\frac{\Gamma\big(\frac{pr}{r-p}\big)}{\Gamma\big(\frac{r}{r-p}\big)\Gamma\big(\frac{p(r-1)}{r-p}\big)}\right)^{\frac{r}{p}-1}\Big(\int_\Omega|\nabla f|^p\,dV\Big)^\frac{r}{p},
\end{eqnarray*}
where $r=\frac{pq}{q-p}$ and $\Gamma(\cdot)$ is the classical gamma function. Simplifying the just-obtained estimates, we get
$$
\Big(\int_\Omega |f|^r\,dV\Big)^\frac1r\le\left(\frac{\Gamma\big(\frac{pr}{r-p}\big)}{\Gamma\big(\frac{r}{r-p}\big)\Gamma\big(\frac{p(r-1)}{r-p}\big)}\right)^{\frac{1}{p}-\frac1r}
\left(\frac{\int_\Omega|\nabla f|^p\,dV}{ E_{n,p}(\Omega)V(\Omega)^{\frac{p}{q}-\frac{p}{n}}}\right)^\frac{1}{p}.
$$
Letting $q\to p$ in the last inequality, we find $r\to\infty$ and thus
$$
\|f\|_{L^\infty(\Omega)}^p\le\big(E_{n,p}(\Omega)V(\Omega)^{\frac{n-p}{n}}\big)^{-1}\int_\Omega|\nabla
f|^pdV,
$$
thereby establishing (iii).

(ii)/(iii)$\Rightarrow$(i) Due to (ii)$\Leftrightarrow$(iii), we may assume that (iii) is valid with $E_{n,p}(\Omega)>0$ (otherwise nothing is to prove). For $f\in C_0^\infty(\Omega)$ and $q>p$ we employ the H\"older inequality to
get
\begin{align*}
\int_\Omega|f|^q\,dV&=\int_\Omega|f|^{q-p}|f|^pdV\\
&\le\left(\Big(\frac{\int_\Omega|\nabla
f|^p\,dV}{E_{n,p}(\Omega)}\Big)^\frac1p
V(\Omega)^{\frac1n-\frac1p}\right)^{q-p}\int_\Omega|f|^pdV
\\[10pt]
&\le\left(\frac{\displaystyle{\int_\Omega|\nabla
f|^pdV}}{\displaystyle{\int_\Omega|f|^pdV}}
\right)^{\frac{q-p}{p}}\big(E_{n,p}(\Omega)\big)^{\frac{p-q}{p}}V(\Omega)^\frac{q-p}{n}\int_\Omega|f|^qdV,
\end{align*}
thereby reaching
$$
\frac{\displaystyle{\int_\Omega|\nabla
f|^pdV}}{\displaystyle{\int_\Omega|f|^pdV}
}\ge E_{n,p}(\Omega)V(\Omega)^{-\frac{p}{n}}.
$$
Furthermore, the formulation of $\lambda_p(\Omega)$ in terms of $C_0^\infty(\Omega)$ (cf. \cite{KaLi}) is used to verify the validity of (i).
\end{proof}

\begin{remark}\label{rm5} The following sharp geometric limit inequality
(cf. \cite[Corollary 15]{KaFr} and
\cite{JuLM}):

\begin{equation*}\label{eqfinal1}
\lim_{p\to\infty}\lambda_p(\Omega)^\frac1p
V(\Omega)^\frac1n\ge\omega_n^\frac1n,
\end{equation*}
along with (\ref{eq10}), induces a purely geometric
quantity
$$
\Lambda_\infty(\Omega):=\lim_{p\to\infty}\gamma_p(\Omega)^\frac1p=\lim_{p\to\infty}\lambda_p(\Omega)^\frac1p=\inf_{x\in\Omega}\hbox{dist}(x,\partial\Omega)^{-1}.
$$
Obviously, (\ref{eq0b}) is used to derive the $\infty$-Faber-Krahn
inequality below:
\begin{equation}\label{eqq}
\Lambda_\infty(\Omega)\ge \Lambda_\infty(\Omega^\ast).
\end{equation}
Moreover, as the limit of $\Delta_p u=\lambda_p(\Omega)|u|^{p-2}u$
on $\Omega$ as $p\to\infty$, the following Euler--Lagrange equation:
$$
\max\{\Lambda_\infty(\Omega)-|\nabla u|u^{-1},\ \Delta_\infty u\}=0~~
\text{in}~~ \Omega
$$
holds in the viscosity sense (cf. \cite{JuLM}), where
$$
\Delta_\infty u(x):=\sum_{j,k=1}^n\Big(\frac{\partial u(x)}{\partial
x_j}\Big)\Big(\frac{\partial^2 u(x)}{\partial x_j\partial
x_k}\Big)\Big(\frac{\partial u(x)}{\partial x_k}\Big)=\langle D^2u(x)\nabla u(x),\nabla u(x)\rangle
$$
is the so-called $\infty$-Laplacian.
\end{remark}

Last but not least, we would like to say that since the geometry of
$\mathbb R^n$ -- the isoperimetric inequality plays a key role in
the previous treatment, the five propositions above may be
generalized to a noncompact complete  Riemannian
manifold (substituted
for $\mathbb R^n$) with nonnegative Ricci curvature and
isoperimetric inequality of Euclidean type, using some methods and
techniques from \cite{Cha,Gri,Heb} and \cite{Sal}.

\vskip7pt
\noindent
\textbf{Acknowledgment.}~The work was partially supported by an NSERC (of Canada) discovery grant and
a start-up fund of MUN's Faculty of Science as well as National Center for Theoretical Sciences (NCTS) located in National Tsing Hua University, Taiwan. The final version of this paper was completed during the author's visit to NCTS at the invitations of Der-Chen Chang (from both Georgetown University, USA and NCTS) and Jing Yu (from NCTS) as well as Chin-Cheng Lin (from National Central University, Taiwan).


\end{document}